\documentclass[12pt]{article}

\usepackage{amsmath,amsthm,amsfonts,amssymb,epsf,epsfig}
\evensidemargin
.25in \textheight 8.5in \topmargin 0in 
\title{Incompressible Maps of Surfaces and Dehn Filling}
\author{ Ulrich Oertel \thanks{Research supported by the National Science Foundation}}

\date{March, 2000; revised January, 2003; preliminary versions 1998}

\newtheorem{thm}{Theorem}[section] 
\newtheorem{lemma}[thm]{Lemma}

\newtheorem{corollary}[thm]{Corollary}
\newtheorem{proposition}[thm]{Proposition}

\theoremstyle{definition}
\newtheorem{defn}[thm]{Definition}
\newtheorem{defns}[thm]{Definitions}

\newtheorem{ex}[thm]{Example}

\theoremstyle{remark}

\begin{document}
\maketitle

\def\reals{\mathbb R}
\def\rationals{\mathbb Q}
\def\complex{\mathbb C}
\def\naturals{\mathbb N}
\def\integers{\mathbb Z}

\def\proj{P}
\def\hyp {\hbox {\rm {H \kern -2.8ex I}\kern 1.15ex}}
\def\intr{\text{int}}
\def\inter{\ \raise4pt\hbox{$^\circ$}\kern -1.6ex}
\def\Cal{\cal}
\def\from{:}
\def\inverse{^{-1}}
\def\Max{\text{Max}}
\def\Min{\text{Min}}
\def\embed{\hookrightarrow}
\def\Genus{\text{Genus}}
\def\Z{Z}
\def\X{X}
\def\roster{\begin{enumerate}}
\def\endroster{\end{enumerate}}
\def\definition{\begin{defn}}
\def\enddefinition{\end{defn}}
\def\subhead{\subsection\{}
\def\theorem{thm}
\def\endsubhead{\}}
\def\head{\section\{}
\def\endhead{\}}
\def\example{\begin{ex}}
\def\endexample{\end{ex}}
\def\ves{\vs}
\def\mZ{{\mathbb Z}}
\def\M{M(\Phi)}
\def\bdry{\partial}
\def\hop{\vskip 0.15in}
\def\trip{\vskip 0.09in}

\centerline{Rutgers University, Newark}

\abstract 
We prove results showing that the existence of essential maps of surfaces in a manifold $\bar M$ obtained from a
3-manifold
$M$ by Dehn filling implies the  existence of essential maps of surfaces in $M$.  
\trip

\noindent {\it Key Words and Phrases:}  Dehn filling, essential surface, incompressible surface, 3-manifold, immersed
surface.
\trip

\noindent {\it Mathematics Subject Classification:}
57N10, 57M99.
\endabstract

\section{Introduction and Statements }\label{Introduction}

The best understood 3-manifolds are Haken 3-manifolds;
these are manifolds that are irreducible and contain incompressible surfaces.
One approach to 3-manifolds is to try to generalize the notion of an incompressible
surface.  Essential laminations and mapped-in $\pi_1$\--injective
surfaces have both been used with some success as substitutes for incompressible
surfaces, yielding theorems which apply to classes of manifolds larger than 
the class of Haken manifolds.

In dealing with maps of surfaces $f\from S\to M$ which induce injections on $\pi_1$, a
proof of the Simple Loop Conjecture would be helpful.  Suppose $M$ is a
compact, orientable closed  3-manifold, $S$ is a closed orientable surface of genus $g\ge
1$, and $f\from S\to M$ is a map with the property that the induced map $\pi_1(S)\to
\pi_1(M)$ has non-trivial kernel.  Then the conjecture is that there is an essential
simple closed curve $\gamma$ in $S$, with $f|_\gamma$ null-homotopic in $M$.  (Note: The same
conjecture with $M$ replaced by a surface was proved by David Gabai in \cite{DG:SimpleLoop}.)

If the conjecture is true, a map  $f\from S\to M$ which is not essential, in the sense
that it does not induce an injection on $\pi_1$, can be replaced by a map of a surface of lower
genus by performing surgery on $S$ and on the map $f$.  In dealing with
embedded incompressible surfaces, this kind of surgery is an important tool.  

For completeness, we give a precise definition of injective and $\bdry$\--injective maps
of surfaces.  We always assume
$M$ is an orientable, compact 3-manifold.  If
$S$ is a compact surface, possibly with boundary, we will say that a map
$f\from (S,\bdry S)
\to (M,\bdry M)$ is {\it $\pi_1$-injective} if the induced map $\pi_1(S)\to \pi_1(M)$ is
injective.  We will say the map
$f\from (S,\bdry S)
\to (M,\bdry M)$ is {\it $\bdry$-$\pi_1$-injective} if the induced map on $\pi_1(S,\bdry S, p)\to \pi_1(M,\bdry M,f(p))$ is
injective for every choice of base point $p$ in $\bdry S$.

One can define a class of maps of
surfaces which would be
$\pi_1$-injective if the Simple Loop Conjecture were true.  
Given an orientable surface $S$ which is not a sphere or disc, we shall say that a
map
$f\from (S,\bdry S)\to (M,\bdry M)$ is {\it incompressible} if no essential simple loop in
$S$ is mapped to a homotopically trivial curve in
$M$.  When $M$ has boundary, we say that the map 
$f\from (S,\bdry S)\to (M,\bdry M)$ is {\it
$\bdry$-incompressible} if no essential simple arc in
$S$ is mapped to an arc in $M$ which is homotopic in $M$, rel boundary, to an arc in $\bdry M$.  If $f$
is a sphere (or a disc) a
map
$f\from (S,\bdry S)\to (M,\bdry M)$ is {\it essential} if it is not null-homotopic (not
homotopic to a disc in $\bdry M$).  We will say a map $f\from (S,\bdry S)\to (M,\bdry M)$
is {\it essential} if either it is an essential map of a disc or sphere, or it is an
incompressible, $\bdry$-incompressible map of a surface which is not a disc or sphere.  We emphasize that
all results in this paper about essential maps of surfaces use the definition of ``essential" in terms of
incompressible maps, i.e. a definition which bypasses the issue of the Simple Loop Conjecture.

It is possible to prove some results about incompressible maps of
surfaces, results analogous to results about embedded incompressible surfaces.  In fact,
in \cite{UO:ImmersedSlopes}, I stated Proposition \ref{OldProp}, below.  Ian Agol pointed out a gap in the
proof.  The purpose of this paper is to generalize the result while filling the gap, a missed case which
appears as Case 5 in the proof of the main theorem of this paper, Theorem \ref{MainThm}.

\begin{proposition}\label{OldProp}  (See \cite{UO:ImmersedSlopes}.) Suppose $M$ is a compact, orientable  3-manifold
with
$k$ torus boundary components, and suppose $\bar M=M(r_1,\ldots, r_k)$ is the
manifold obtained by performing Dehn filling of slope $r_i$ on the i-th boundary
torus.  Suppose
$S$ is a closed orientable surface, not a sphere.  Given an incompressible (e.g. $\pi_1$-injective)
map $f\from S \to\bar M$, there is an incompressible,
$\bdry$-incompressible map $\hat f\from (\hat S,\bdry \hat S) \to (M,\bdry M)$,
where $\hat S$ is a surface obtained from
$S$ by removing interiors of some collection of disjointly embedded discs and where
components of
$\bdry
\hat S$ mapped to the $i$-th boundary component are embedded curves of slope $r_i$.  
\end{proposition}

Agol asked me about the following related result,
which he believed to be true.  He proved a somewhat weaker version, Lemma 3.2 in \cite{IA:Bounds}.  

\begin{proposition} \label{AgolProp}  Suppose $M$ is an orientable 3-manifold with $\bdry
M$ a single torus, let $M(r)$ be the manifold obtained from $M$ by performing a
slope $r$ Dehn filling on $M$, and let $K$ be the core curve of the Dehn filling solid
torus $R$.  Suppose
$f\from S\to M(r)$ is a map of a disc to $M(r)$ with $f|_{\bdry S}$ mapped to $K$ with
degree $\ge 1$. 
Choosing a map which intersects $K$ minimally yields an incompressible,
$\bdry$-incompressible map
$\hat f\from
(\hat S,\bdry\hat S)\to (M,\bdry M)$ of a punctured disc $\hat S$, where $\hat
S=f\inverse(M)$.  One boundary component of $\hat S$ is mapped to a curve in $\bdry M$ homotopic
in $R$ to $f|_{\bdry S}$, and all other curves of $\bdry \hat S$ are embedded as curves of slope $r$.
\end{proposition}

Exactly what we mean by a map ``which intersects $K$ minimally" will be explained below.
For Agol, the motivation for proving the above result was to be able to conclude that if 
a Dehn filling on $M$ yields a finite fundamental group 3-manifold, then there exists an
incompressible,
$\bdry$-incompressible map of a disc-with-holes into the manifold $M$.

One can prove a more general theorem which has all of the above as corollaries.
To state the more general theorem, we first make some more definitions. Let $K=\cup K_j$ be
a link in a manifold $\bar M$, where each $K_j$ is a knot, and let $R_j$ be a solid torus
regular neighborhood of $K_j$. Suppose
$S$ is a surface and suppose $\bdry S$ is the disjoint union of $\bdry_i S$ ({\it
interior boundary}) and
$\bdry_eS$ ({\it exterior boundary}).

\begin{defns} A map $f\from (S,\bdry_iS, \bdry_eS)\to (\bar
M,K,\bdry \bar M)$ is {\it relatively incompressible} at
$K$, if 

i) for no simple essential arc $\alpha$ in $S$ with both ends on $\bdry_i S$ is $f|_\alpha$ homotopic in
$\bar M$ into $K$.

ii) for every component $\alpha$ of
$\bdry_iS$
$f|_\alpha$ is not null-homotopic in $K$, i.e. has $\text{degree}\ge 1$.  

The map
$f\from (S,\bdry_iS, \bdry_eS)\to (\bar
M,K,\bdry \bar M)$ is {\it incompressible} if (as before) no essential simple loop in
$S$ is mapped to a homotopically trivial curve in
$\bar M$.  When $\bar M$ has boundary, we say that the map 
$f\from (S,\bdry_iS, \bdry_eS)\to (\bar
M,K,\bdry \bar M)$ is {\it
$\bdry$-incompressible} if no essential simple arc $(\alpha,\bdry \alpha)\in (S,\bdry_e S)$ is mapped to an
arc in $\bar M$ which is homotopic in $\bar M$ to an arc in $\bdry \bar M$.  In all of the following, it is
helpful to think of the special case where $\bar M$ is closed and $\bdry_eS=\emptyset$.

As we shall see, it is possible to homotope a relatively incompressible map
$f\from (S,\bdry_iS, \bdry_eS)\to (\bar
M,K,\bdry \bar M)$ such that
$f\inverse(\cup R_j)$ is a collar of $\bdry_iS$,
together with discs mapped homeomorphically to meridians of $R_j$'s.  In this case we say $f$ is a {\it good
map} with respect to $K$.  For a good map, we also require that $f|_{\bdry_eS}$ be in general position,
i.e. that there be no triple points on $\bdry \bar M$.  Similarly, we require $f$ to intersect $\cup_j \bdry
R_j$ generically.
\end{defns}
\trip

To see that a relatively incompressible map $f\from (S,\bdry_iS, \bdry_eS)\to (\bar M,K,\bdry
\bar M)$ can be replaced by a homotopic good map, we proceed as follows.
Choose a curve
$\gamma$ of $\bdry_iS$ and suppose it is mapped to $K_j$.  Choose any curve
$\delta$ embedded in
$\bdry R_j$ homotopic in $R_j$ to the map $f|_\gamma$ and intersecting meridian discs efficiently. We can then
construct a map of an annulus
$A$ into $R_j$ by coning to the center
$\delta\cap E$ in each meridian disc $E$ of a 1-parameter family of meridian discs which
foliate $R_j$.  Then we homotope $f$ so that on a collar of
$\gamma$ it agrees with the map of $A$ just constructed.  We do this for every component of
$\bdry_iS$.  Finally, making
$f$ transverse to $\cup_jK_j$ on $S-N(\bdry_iS)$, we may assume that remaining non-annular
components of
$f\inverse(\cup_jR_j)$ are discs mapped to meridians of $R_j$'s.  This yields a good
map after a further homotopy to achieve general position in $\bdry \bar M$ and in $\cup_j\bdry R_j$.

We will use $M$ to denote $\bar M-\intr(\cup_j R_j)$, $\hat S=f\inverse(M)$, and $\hat f=f|_{\hat S}$. 

It will be convenient to extend our definition of essentiality.  

\begin{defn}A map 
$$f\from (S,\bdry_iS, \bdry_eS)\to (\bar M,K,\bdry \bar M)$$ of an orientable surface
(possibly a disc or sphere) is {\it essential} if one of the following holds:

\begin{narrower}

1) It is not a map of a disc or sphere and is incompressible, $\bdry$-incompressible, and
relatively incompressible at $\cup K_j$.

2) It is an essential map of a sphere or an essential map of a disc $(S,\bdry_eS)\to (\bar M,\bdry\bar M)$.

3) It is a relatively incompressible map of a disc $f\from (S,\bdry_iS)\to (\bar M, \cup_j K_j)$, i.e., a map
of a disc $S$ of non-zero degree on $\bdry S$.  

\end{narrower}
\end{defn}

We assign a complexity to a good map $f\from (S,\bdry_iS,\bdry_eS)\to$ $ (\bar
M,\cup_jK_j,\bdry\bar M)$  The complexity is $c(f)=m(f)+i(f)+e(f)$, where $m(f)$ is the
number of discs in $f\inverse(\cup_jR_j)$, $i(f)$ is the number of
double points of $f$ on $\cup_j\bdry R_j$, and $e(f)$ is the number of double points
on $\bdry \bar M=\bdry_eM$.

\begin{thm}\label{MainThm} Suppose $M$ is an orientable 3-manifold with
$k$ distinguished torus boundary components $T_j$, and suppose $M(r_1,\ldots, r_k)=\bar M$ is the
manifold obtained by performing Dehn filling of slope $r_j$ on the $j$-th boundary
torus $T_j$.  Let $K_j$ denote the core curve of the $j$-th Dehn filling solid torus,
$R_j$.    Suppose
$S$ is an orientable, connected surface and suppose $f\from (S,\bdry_i S, \bdry_e S)\to
(\bar M,\cup_j K_j,
\bdry \bar M)$, is essential.  Then if $c(f)$ is minimal over all
essential good maps  
$g\from (S,\bdry_i S,
\bdry_e S)\to (\bar M,\cup_j K_j,
\bdry \bar M)$ with $g$ homotopic to $f$ on $\bdry_iS$, then the restriction $\hat f$ of $f$ to
$\hat S=f\inverse(M)$ is essential in $M$.  For every component of $\bdry_i S$ mapped to $K_j$ there is a
curve of
$\bdry
\hat S$ isotopic to it in $S$ which is mapped by $\hat f$ to $\bdry R_j$.  Every other boundary curve of
$\hat S$ in
$f\inverse (\bdry R_j)$ is embedded by $\hat f$ as a meridian in $\bdry R_j$. 
\end{thm}

Note that, in general, the theorem yields an essential
map $\hat f$ of a surface $\hat S$ whose boundary can be subdivided into three parts.
The exterior boundary of $S$ remains unchanged in $\hat S$, yielding $\bdry_e\hat S$.  The interior boundary
of
$\hat S$,  $\bdry_i\hat S$, consists of those closed curves of $\bdry \hat S$ on the boundary of annular collars in
$S-\hat S$ of
$\bdry S$ mapped by $f$ to
$\cup_j R_j$.  The components of $\bdry_i \hat S$ are isotopic in $S$ to components of $\bdry_i S$.  Finally, $\bdry_m
\hat S=\cup f\inverse(\bdry R_j)-\bdry_i\hat S$.  In many applications, we can take
$\bdry_eS=\emptyset$. In case $\hat f$ is an embedding on a component of $\bdry_i \hat S$,
we can assign a {\it boundary slope} to that boundary component.  The terminology is
standard:  A slope is the isotopy class of an embedded closed curve in a torus.  In Proposition \ref{AgolProp} above
we obtain two boundary curves, one embedded and yielding a boundary slope.  In Corollaries \ref{HomologyCor} and
\ref{KnotCor} below, we obtain exactly two boundary slopes, one of which is the slope of the Dehn filling. This
indicates that incompressible,
$\bdry$-incompressible maps of connected surfaces with multiple boundary slopes, as
described in  \cite{UO:ImmersedSlopes}, occur quite naturally.

\begin{corollary} \label{HomologyCor} Suppose $M$ is an orientable 3-manifold (for example the exterior
of a knot in $S^3$) with
$\bdry M$ a single torus.  Let $M(r)$ be the manifold obtained from $M$ by performing a
slope $r$ Dehn filling on $M$, and let $K$ be the core curve of the Dehn filling solid
torus.  If $K$ is null-homologous in $M(r)$, then $K$ bounds a mapped-in surface $f\from S\to M(r)$ (so $f|_{\bdry
S}$ is an embedding to $K$).  If
$S$ has genus $g$, then there is an essential map
$\hat f\from
\hat S\to M$, of a surface $\hat S$ of genus $g$ , with
$\hat f$ mapping one component, $\bdry_i\hat S$,  of $\bdry \hat S$ to an integer slope embedded simple closed
curve in
$\bdry M$ and embedding the remainder, $\bdry_m\hat S$, of
$\bdry\hat S $ as slope $r$ closed curves in
$\bdry M$.

In particular, if $M(r)$ is a homology sphere, there is a map $\hat f$ of a surface $\hat
S$ as above.
\end{corollary}

\begin{corollary}\label{KnotCor}  Let $K\subset S^3$ be a knot in the 3-sphere and let $X$ be
the knot exterior.  There is an essential map
$\hat f\from (P,\bdry P)\to (X,\bdry X)$ of a planar surface, with
$\hat f$ restricted to one component of $\bdry P$, $\bdry_iP$ say, mapping to an integer
slope embedded simple closed curve in $\bdry X$ and with all other components of $\bdry P$
mapping to embedded meridian curves in
$\bdry X$.
\end{corollary}

More general than Corollary \ref{HomologyCor}, we have the following.

\begin{corollary} \label{HomologyCorTwo} Suppose $M$ is an orientable 3-manifold with
$k$ torus boundary components, and suppose $\bar M=M(r_1,\ldots, r_k)$ is the
manifold obtained by performing Dehn filling of slope $r_j$ on the j-th boundary
torus $T_j$, for each $j$.  Let $K_j\subset \bar M$ denote the core curve of the $j$-th Dehn filling solid
torus. Giving an orientation to each $K_j$, suppose $\sum_i n_{j_i}[K_{j_i}]=0$ in
$H_1(\bar M)$ for some integers
$n_{j_i}$, and no proper subset of the $[K_{j_i}]$'s is dependent in
$H_1(\bar M)$.  Then the 1-cycle $\sum_i n_{j_i}[K_{j_i}]$ bounds a connected mapped-in surface $f\from S\to
\bar M$.  If the genus of $S$ is 
$g$, then there is an essential map
$\hat f\from
\hat S\to M$ of a connected surface $\hat S$ of genus $g$ with the following property.  The boundary of $\hat S$
can be decomposed as $\bdry
\hat S=\bdry_i
\hat S\cup\bdry_m
\hat S$, so that for every $T_{j_i}$ there is at least one component of $\bdry_i\hat S$ mapped to $T_{j_i}$,
and so that every component of $\bdry_m\hat S$ is embedded in some $T_j$ as a slope $r_j$ curve.
\end{corollary}

The above corollary is related to Lemmas 2.4 and 2.5 of \cite{ML:HyperbolicDehn} which deal with planar surfaces
only.

I thank Ian Agol, with whom I exchanged ideas related to this paper.   The main result
presented here grew out of Proposition \ref{AgolProp}, which was formulated by Agol.

\section{Proof of Main Theorem}\label{MainProof}

\begin{proof} {\it(Theorem \ref{MainThm}.)}  Recall that $R_j$ denotes the $j$-th filling solid torus attached
to the 2-torus $T_j\subset \bdry M$.  

We assume now that $f$ is chosen so that the complexity $c(f)$ is minimal over all good
maps
$g\from (S,\bdry_i S,
\bdry_e S)\to (\bar M,\cup_j K_j,
\bdry \bar M)$ with $g$ homotopic to $f$ on $\bdry_iS$.  Henceforth we let
$\hat f$ denote $f$ restricted to
$\hat S=f\inverse(M)$.  We shall deal with the case that $\hat f\from \hat S\to M$ is a
map of a disc or sphere at the end of the proof, showing that the map is essential.  In 
all remaining cases, our task is to prove incompressibility and $\bdry$-incompressibility of
the map $\hat f$.

If $\hat f$ is compressible, then using the incompressibility
of $f$, we can replace $f$ by a map whose image includes fewer meridian discs in
$R_j$'s as follows.  A simple closed curve $\alpha\embed \hat S$ with $\hat f|_\alpha$
null-homotopic in $M$, but with $\alpha$ not null-homotopic in $\hat S$, must be
null-homotopic in $S$ by the incompressibility of $f$.  Replacing $f$ restricted to the
disc $D$ in $S$ bounded by $\alpha$ by the null-homotopy in $M$, we obtain $f$
intersecting $\cup R_j$ in fewer discs, so the summand $m(f)$ of $c(f)$ is reduced, a contradiction.
The same argument applies if $S$ is a disc or a sphere.

Showing that $\hat f$ is $\bdry$-incompressible is considerably more
difficult.  Suppose $\beta$ is a simple essential arc in $\hat S$ joining points in 
$\bdry\hat S$ and homotopic to an arc in $\bdry M$.  The homotopy is represented by
a map $h\from (H,a) \to (M,\bdry M)$, where $H$ is a half-disc, i.e. a disc with $\bdry H$ divided into two
complementary arcs $a$ and $b$ intersecting only at their endpoints.  The map 
$h|_b$ factors through an embedding
$k\from b\embed \hat S$ whose image is
$\beta\subset \hat S$, so
$\hat f\circ k=h|_b$.  We may suppose the arc $h|_a$ is not homotopic in $\bdry M$ to an arc in the
image of
$\hat S$, or more precisely, that $h|_a$ is not homotopic to $\hat f$ restricted to an arc in $\bdry
\hat S$, otherwise extending this homotopy to $h$ changes $h\from H\to M$ into a potential compressing disc for
$\hat f$.  The fact that $\hat f$ is incompressible then shows that $h$ is not a $\bdry$-compression, since
the image $\beta$ of $h|_b$ is then inessential in $\hat S$.

Keeping in mind that we do not need to consider the uninteresting $\bdry$- compressions described above, there
are 6 cases.  We shall describe the cases by homotoping $\hat f|_{\beta}$ in $M$ to an arc near $\bdry M$,
extending the homotopy to all of $S$ rel $\bar M -\inter M$ to get a homotopy of $f$, and homotoping $f$ further to
push
$\beta$ into an
$R_j$.  Case 1 and 2, shown in Figures \ref{DehnFold} and \ref{DehnTwistFold}a, occur when the arc
$\beta$ joins boundary components of
$\hat S$ which are capped by meridian discs in $S$.  Each figure depicts a portion of the image of
$f$ in an $R_j$, or in a regular neighborhood of an $R_j$.  Case 3, shown in Figure \ref{DehnCollarJoin},
occurs when $\beta\subset\hat S$ joins a collar of $\bdry_iS$ to itself, or joins a pair
of annular collars.  Case 4, shown in Figure \ref{DehnBetaToCollar}a, occurs when $\beta$ joins an annular
collar
to a meridian disc.  Case 5, shown in Figure \ref{DehnLongBand}, occurs when $\beta$ joins a meridian disc
to itself.  Case 6 is the
less interesting case in which $\beta$ joins $\bdry_e\hat S$ to itself and $\beta$ is
homotopic in $M$ to an arc in $\bdry_eM=\bdry\bar M$.  Note that all of the figures show a
regular neighborhood in $\hat S$ of
$\beta$ together with the collars or meridian discs they join; we do not show all of
$f(S)\cap (\cup R_j)$.

Our list of cases is exhaustive:  Since both ends of $\beta$ must be mapped to the
same component of $\bdry M$, either both ends lie $\bdry_eM$, or both ends lie in a torus
component $\bdry R_j$ of
$\bdry_iM=\bdry M-\bdry_eM$.  If both ends lie on $\bdry_eM$, we have Case 6. 
Otherwise each end lies either on the boundary of a meridian disc or on the boundary
of a collar in some $R_j$. Cases 1 and 2 deal with arcs $\beta$ joining a meridian disc
to a different meridian disc.  Case 5 deals with arcs $\beta$ joining a meridian disc
to itself.  Case 3 deals with an arc $\beta$ joining two different collars, or a collar
to itself. Case 4 deals with an arc $\beta$ joining a meridian disc to a collar.  One
might think that it is necessary to consider cases like those shown in Figures \ref{DehnFold},\ref{DehnCollarJoin},
or \ref{DehnBetaToCollar} where connecting bands have different numbers of half-twists, but homotopy of
$f$ allows the elimination of pairs of half-twists (of the same sense).  In Figure \ref{DehnBetaToCollar}a,
introduction of a half-twist in the connecting band yields a band which can be homotoped
to an untwisted band on the opposite side of the arc of intersection.

\begin{figure}[ht]
\centering
\scalebox{1.0}{\includegraphics{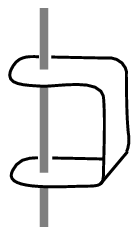}}
\caption{\small Case1.}
\label{DehnFold}
\end{figure}

\begin{figure}[ht]
\centering
\scalebox{1.0}{\includegraphics{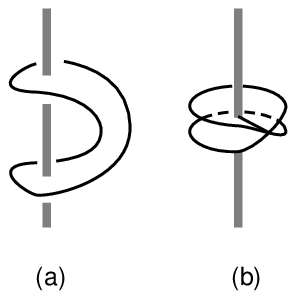}}
\caption{\small Case 2. }
\label{DehnTwistFold}
\end{figure}

\begin{figure}[ht]
\centering
\scalebox{1.0}{\includegraphics{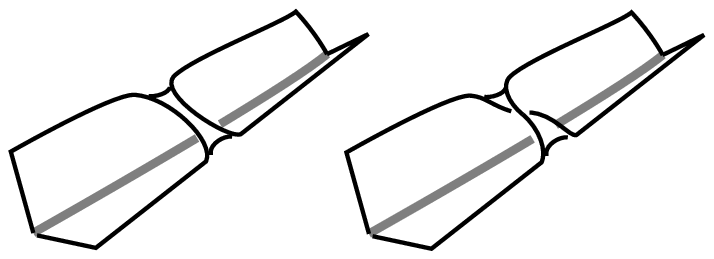}}
\caption{\small Case 3.}
\label{DehnCollarJoin}
\end{figure}

\begin{figure}[ht]
\centering
\scalebox{1.0}{\includegraphics{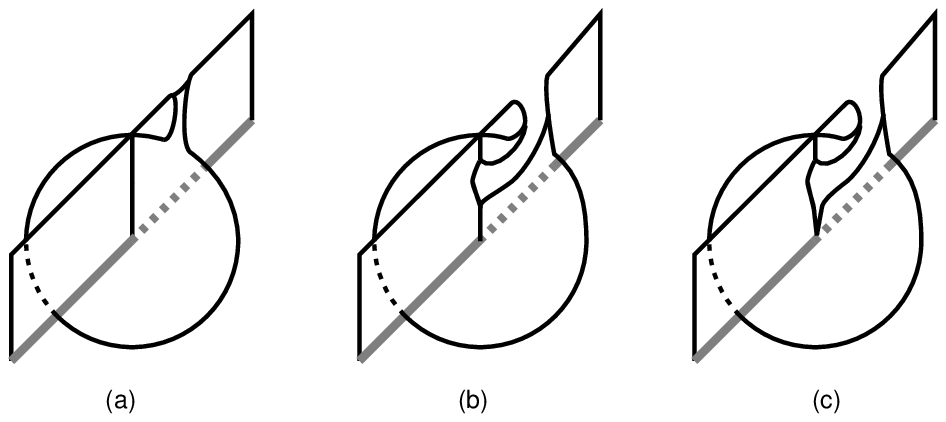}}
\caption{\small Case 4.}
\label{DehnBetaToCollar}
\end{figure}

\begin{figure}[ht]
\centering
\scalebox{1.0}{\includegraphics{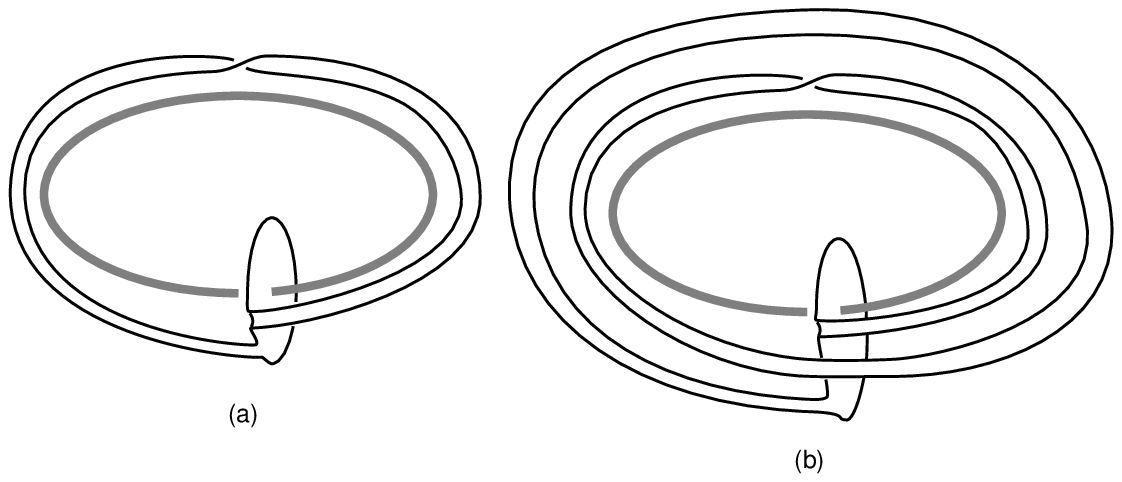}}
\caption{\small Case 5.}
\label{DehnLongBand}
\end{figure}

We will obtain a contradiction in each case.
\bigskip

\noindent {\it Case 1.} The arc $\beta$ joins two meridian discs of $S$ with opposite
orientations, Figure \ref{DehnFold}. In this case, $f$ can
be homotoped to reduce the number of disc intersections of $f(S)$ with $\cup_i R_i$, a
contradiction. 

\bigskip
\noindent {\it Case 2.} The arc $\beta$ joins two meridian discs $E_1$ and $E_2$ of
$f(S)\cap R_j$, for some $j$, with consistent orientations, Figure \ref{DehnTwistFold}a.  
In this case, the map $f$ can be homotoped such
that the image of a neighborhood of 
$E_1\cup E_2\cup \beta$ is as shown in Figure  \ref{DehnTwistFold}b.  The homotopy can be done in
such a way that
$f(S)$ is otherwise unchanged in a neighborhood of $\cup_i R_i$.  So we have
replaced two discs of $f\inverse(\cup_i R_i)$ by one disc mapped with an order 2 branch point, which we
label
$P$.  We can move the branch point slightly away from $R_j$, if we wish, so that we can still regard $f$ as
having the same complexity $c(f)$.

In order to obtain a contradiction, we must do some global modifications of the map
$f$.  

Applying H. Whitney's Theorem 2 of \cite{HW:DifferentiableManifolds} as described in Lemma 1 of \cite{HW:Immersions2},
we can homotope $f$ such that the
image is locally modelled as shown in Figure \ref{DehnModels}.  At each point of the image, it is
locally modelled as an embedding, as a transverse intersection along an arc, as an order 2 branch point, which
is the cone on a figure eight, or as a triple point.  We may assume that the surface $S$ is transversely
oriented; in the figures we draw intersections so that the angle from the ``+-side" of one
sheet of $S$ to the +-side of the other sheet of $S$ is obtuse.  Throughout this paper a {\it branch point}
will be a branch point of order 2.

\begin{figure}[ht]
\centering
\scalebox{1.0}{\includegraphics{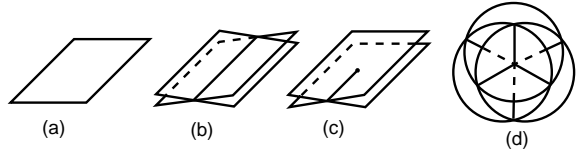}}
\caption{\small Local models for surface image of surface in $M$.}
\label{DehnModels}
\end{figure}

The idea is to eliminate branch points using
homotopies starting at a branch point.  Given an arc of self-intersection
with at least one end a branch point, there is an obvious homotopy which moves
the branch point along the arc of intersection while shortening the arc, see Figure \ref{DehnZip}.  The
effect is similar to the effect of a zipper, but four rather than two sheets of cloth come
together along the zipper or arc of intersection. The usual zipper cuts a surface or rejoins it;
our zipper cuts two surfaces and rejoins them locally in the opposite way.   The
branch point corresponds to the zipper head.  The zipping operation was used by H. Whitney to prove Theorem 6 of 
\cite{HW:Immersions2} and is also
used in Bing's monograph, \cite{RHB:ThreeManifolds}. In the proof of the Simple Loop Conjecture for maps from surfaces
to surfaces, see \cite{DG:SimpleLoop}, and in \cite{DGWK:SurfaceMaps} the zipping operation was referred to as ``the
calculus of double curves." 

\begin{figure}[ht]
\centering
\scalebox{1.0}{\includegraphics{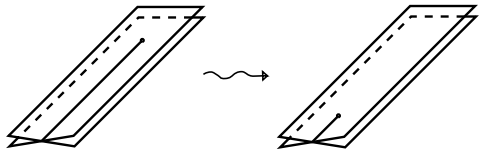}}
\caption{\small The zipping homotopy.}
\label{DehnZip}
\end{figure}

\begin{figure}[ht]
\centering
\scalebox{1.0}{\includegraphics{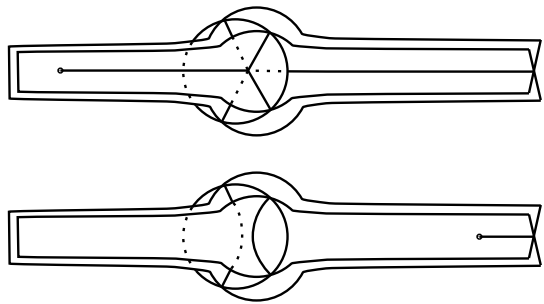}}
\caption{\small Zipping through a triple point.}
\label{DehnZipTriple}
\end{figure}

We begin by performing zipping operations with the goal of eliminating the branch point $P$.  We choose any {\it other} branch point and begin zipping. It 
is possible that $P$ is the only branch point, but we shall deal with that eventuality
later.
If a triple point is encountered, we continue to zip through the triple
point as shown in Figure \ref{DehnZipTriple}.  Clearly the curves of self-intersection of $S$ at
the triple point are altered when the branch point passes through the triple
point; there is a surgery operation on the double curves.  Notice that in our
figures we use a standard model for a triple point, as shown in Figure \ref{DehnModels}d. 
This explains our way of drawing Figure \ref{DehnZipTriple}.  After passing through a triple
point, the zipper may return to what was the triple point once more, or the
zipping path could return to the same point a second time as well.

There are just four possible obstructions which make it impossible to continue
zipping: namely,  the zipping arc may end at a point on some $K_j$, where a meridian of
$f(S)\cap R_j$ meets the image of $\bdry_i S$, as shown in Figure \ref{DehnZipCore}; it may end at
another branch point as shown in Figure \ref{DehnZipSurgery}; it may end at $\bdry \bar M$, where
$f(\bdry S)$ has a self-intersection as shown in Figure \ref{DehnZipToBdry}; or it may end where two points
of $\bdry_iS$ are mapped to the same point of a $K_j$ at the end of an arc of intersection,
as shown in Figure \ref{DehnZipToK}.

\begin{figure}[ht]
\centering
\scalebox{1.0}{\includegraphics{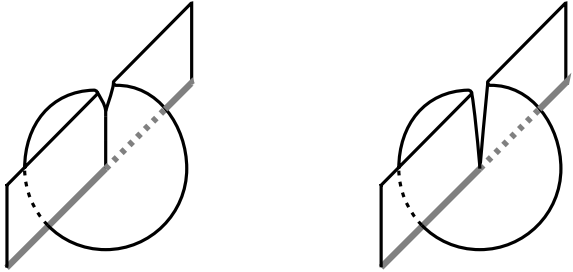}}
\caption{\small Zipping ending at a meridian in $R_j$.}
\label{DehnZipCore}
\end{figure}

\begin{figure}[ht]
\centering
\scalebox{1.0}{\includegraphics{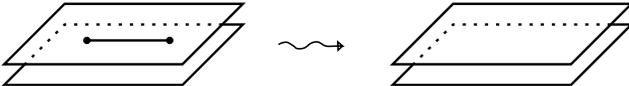}}
\caption{\small Zipping ending at branch point.}
\label{DehnZipSurgery}
\end{figure}

\begin{figure}[ht]
\centering
\scalebox{1.0}{\includegraphics{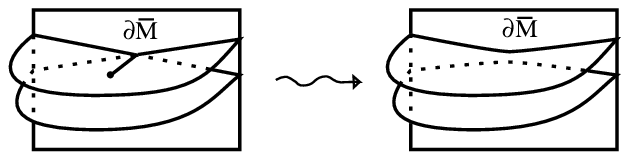}}
\caption{\small Zipping ending at $\bdry\bar M$.}
\label{DehnZipToBdry}
\end{figure}

\begin{figure}[ht]
\centering
\scalebox{1.0}{\includegraphics{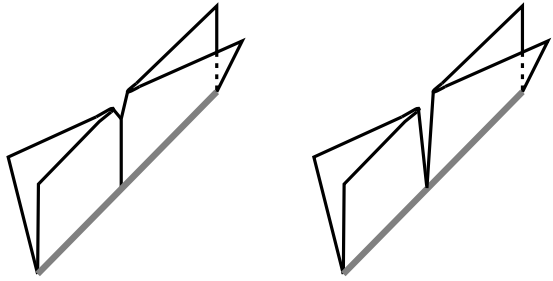}}
\caption{\small Zipping ending at collar intersection at $K$.}
\label{DehnZipToK}
\end{figure}

In the first case, when the zipping ends at a $K_j$ traversing an arc of intersection
between a meridian disc and a collar in an $R_j$, Figure \ref{DehnZipCore} shows that one
obtains a new map of the same surface intersecting
$R_j$ in fewer meridian discs.  The new map, since it is homotopic to the old one and is good, satisfies the
hypotheses of the theorem, which contradicts our assumption that
$f$ has minimal complexity, since the summand $m(f)$ of the complexity $c(f)$ is reduced.

In the second case, when the zipping ends at another branch point, we cut on the arc
of intersection and paste to obtain a surface which is embedded in a neighborhood of the
arc, see Figure \ref{DehnZipSurgery}.  The final surgery respects orientation.  In terms of the topology of
the surface $S$, this final cancellation of branch points locally has the effect of
replacing an annulus of the source surface by two discs; thus it performs a
surgery.   But in our situation, the core of the annulus is a simple closed curve
embedded in $S$ hence it must bound a disc in $S$ by the
incompressibility of the map $f$.  It follows that the final cancellation
splits off a sphere, which we discard if $S$ is not a disc or sphere.  This does not increase the number of
disc components of $f\inverse(\cup R_j)$, and it may reduce it. We obtain a new map $f$
defined on a homeomorphic surface $S$ with fewer branch points.  The new map satisfies the hypotheses of the theorem.
We shall see that this is the only case which does not immediately contradict the minimality of the complexity
$c(f)$.  In effect, we are using a complexity $(c(f),b(f))$ with lexicographical order, where $b(f)$ is the number of
branch points, and in this case $b(f)$ is reduced. 

If the
zipping ends at
$P$, then we have clearly removed the branch point, and we have reduced the number $m(f)$ of 
disc components of $f\inverse(\cup R_j)$ by at least one, as compared with the number of
intersections before the introduction of the branch point $P$.  Again, the new map satisfies the hypotheses of the
theorem, and we have a contradiction to the minimality of the complexity $c(f)$.  If
$f$ is a map of a disc or sphere, the surgery resulting from the cancellation of singularities yields maps of two
spheres, or of a disc and a sphere.  In each case, one of the new maps must be essential on a homeomorphic surface
(disc or sphere), and it either has smaller complexity $c(f)$ or the same $c(f)$ and smaller $b(f)$,
with strictly smaller $c(f)$ if the cancellation involves the branch point $P$.

In the third case, the zipping ends at $\bdry \bar M$, see Figure \ref{DehnZipToBdry}.  Zipping to the end
of the arc of intersection has the effect of surgering $S$ on a simple properly embedded
arc. The $\bdry$-incompressibility of $f$ shows that the surgery cuts a disc from $S$,
so we obtain a new map on a homeomorphic surface by restricting to the remainder of the
surface.  The new map satisfies the hypotheses of the theorem.  This is clear except in the case where
$S$ is a disc, in which case the surgery yields two maps of discs, at least one of which is essential.  In all cases
we have reduced the summand
$e(f)$ of the complexity
$c(f)$, a contradiction.  

In the fourth case, the zipping ends at a $K_j$, see Figure \ref{DehnZipToK}, after traversing an arc of
intersection between two collars in an $R_j$.  Zipping to the end of the arc of
intersection has the effect of surgering
$S$ on a simple properly embedded arc. The  relative incompressibility of $f$ shows that
the surgery cuts a half-disc from
$S$, so we obtain a new map on a homeomorphic surface, satisfying the hypotheses of the theorem, by restricting to the
remainder of the surface.  We have reduced the summand $i(f)$ of the complexity $c(f)$, a contradiction.  As before,
if $S$ is a disc, one of the discs resulting from the surgery will be essential.

We continue zipping
operations starting at branch points other than
$P$.  These operations must all end with cancellations of pairs of branch points.  
The fact that each zipping operation reduces the number $b(f)$ of branch points ensures that the process must stop
when image of $f$ has no branch points or just the branch point at $P$.

If after a sequence of zipping operations, $P$ is the only remaining branch point, then we
perform a zipping starting at $P$.  First this has the effect of replacing the disc
component of $f\inverse(\cup R_j)$ containing $f\inverse(P)$ by two discs again, 
returning the configuration near $P$ to what it was before the introduction of the
branch point 
$P$.  The zipping now must end at a $K_j$ as shown in Figures 9 and 12, or it must end at
$\bdry
\bar M$ as shown in Figure 11.  In each case, one obtains a surface satisfying the hypotheses of the theorem
with reduced complexity as before, contradicting our assumption that $c(f)$ is minimal.

\bigskip
\noindent {\it Case 3.}  The arc $\beta$ joins two annular collars adjacent to components
of $\bdry_iS$ or joins an annular collar to itself.  Extending
$\beta$ yields an arc $\beta'$ in $S$ joining components of $\bdry_iS$, which is
homotopic into
$K_j$.  Relative incompressibility implies that $\beta'$ is not essential in $S$; it cuts
a half-disc $H$ from $S$.  This yields a punctured half-disc in $\hat S$ and shows that $\beta$ joins an
annular collar to itself.  If there is at least one puncture in $H$, then replacing $f$ on the half-disc by the
homotopy from
$\beta'$ to $K_j$ shows that $f$ does not intersect $\cup K_j$ minimally, showing that $c(f)$ was
not minimal.  Thus
$H$ is embedded in $\hat S$, and $\beta$ is not essential in $\hat S$, contrary to
assumption.

\bigskip
\noindent {\it Case 4.}  The arc $\beta$ joins an annular collar adjacent to a component
of $\bdry_iS$ to a meridian disc.  In this case, a neighborhood in $\hat S$ of $\beta$ together with the
meridian disc and the annular collar are as shown in Figure \ref{DehnBetaToCollar}a.  The map $f$ can then be
homotoped to eliminate a meridian disc as shown in Figure \ref{DehnBetaToCollar}b,c.  First, a homotopy introduces a
pair of branch points, then one of the branch points is zipped to $K$.  This again violates our assumption that $f$
has minimal complexity
$c(f)$.

\bigskip
\begin{figure}[ht]
\centering
\scalebox{1.0}{\includegraphics{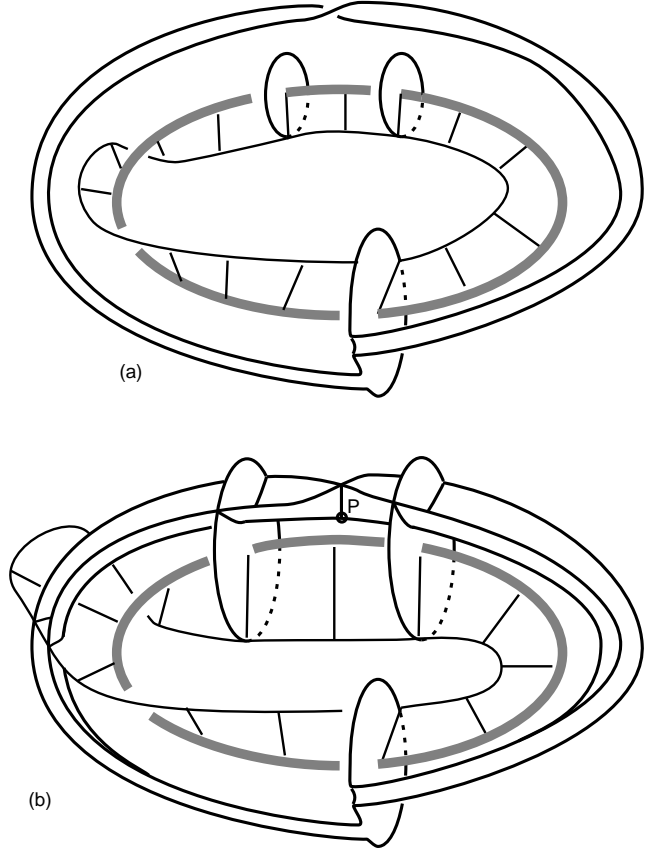}}
\caption{\small Standard $f(A)$.}
\label{DehnLongWithSingBand}
\end{figure}

\bigskip

\noindent {\it Case 5.}  We have a disc $E$ of $f(S)\cap R_j$ joined
to itself by an arc $\beta$ embedded in $\hat S$, and $f(\beta)$ can be homotoped to an 
arc $\beta'$ in $\bdry R_j$ as shown in Figure \ref{DehnLongBand}.  The arc $\beta'$ may wrap around $R_j$ 
longitudinally any number of times.  We need not consider the possibility that $\beta'$ is homotopic to an
arc in $\bdry E\subset \bdry \hat S$, since, as we have observed, in this case the $\bdry$-compression
yields a compressing disc.  The arc $\beta'$ should be made simple at the expense of introducing additional
intersections with $\bdry E\subset \bdry \hat S$.  This can be done by pushing self-intersections of $\beta'$ off
the ends of $\beta'$.  Letting $N(R_j)$ be a regular neighborhood of $R_j$, in effect a slightly larger $R_j$,
we will homotope
$f$ so that
$f(S)\cap N(R_j)$  contains a neighborhood of $f(N(\beta))$, with the image in a kind of general position.  After
the homotopy, $\beta$ is the same as $\beta'$, but pushed out of $R_j$ slightly.  Here
$N(\beta)$ is a neighborhood of $\beta$ in $\hat S$ which we shall refer to as the ``band" connecting $E$ to
itself.  The image
$f(S)\cap N(R_j)$ will consist of the union of some  meridian discs $E_i$, 
the images of annular collars of components of $\bdry_i S$ mapped to $N(R_j)$ , together with
$E$ and the band connecting
$E$ to  itself.  We are assuming that for every disc $E_i$ of $f(S)\cap R_j$ there is a corresponding disc in
$N(R_j)$, and, abusing notation, we denote this larger disc by the same symbol.  The disc
$E$ with a band connecting
$E$ to itself is the image of an annulus 
$A$ in $S$.  We homotope $f$ such that $f$ is a proper map of $(A,\bdry A)$ to
$(N(R_j),\bdry N(R_j))$, and such that the image has exactly one branch point, which we
denote $P$, as shown in Figure \ref{DehnLongWithSingBand}b.  The branch point arises because the band joining
$E$ to itself must have a half-twist. Figure \ref{DehnLongWithSingBand} shows how to make $f(A)$ standard. 
In the figure, the image of the annulus $A$ intersects the meridian discs $E_i$ only in
embedded arcs, and it has one similar embedded arc of self-intersection for every
intersection of the interior of
$\beta$ with
$\bdry E$.  
In addition, there is one arc of self-intersection of
$f(A)$ with one end at
$P$ and the other end at a self-intersection of $f(\bdry A)\subset\bdry N(R_j)$.  Finally,
there are simple arcs of intersection where $A$ intersects annular collars.  These arcs are essential in
$A$, but inessential in the collars.  (The figure shows just one of these arcs, but there may be many.)  The
band may not wrap around $\bdry R_j$ ``monotonically in the longitudinal direction" as shown in the figures,
because we homotoped $\beta'$ to be embedded.  If
$f(\beta')$ intersects $E$ $n$ times, then the band joining $E$ to itself intersects $E$ in $n-2$ arcs.  We use
$n=n(\beta')=n(\beta)$ as a complexity for a $\bdry$-compression of the type we are considering in this case. 

Once again, we homotope $f$, relative $N(R_j)$, such that the image is locally isomorphic
to one of the models of Figure \ref{DehnModels} in the remainder of $\bar M$.  As
before, we begin zipping at a branch point other than $P$, if possible. 

In the course of zipping, it may happen that, in
$N(R_j)$, some arc of intersection disjoint from $\cup K_j$ (other than the arc including $P$) is eliminated by
the zipping.  The zipping effects an orientation respecting cut-and-paste on the arc.  If
this happens, we immediately stop
zipping.  The eliminated arc of intersection was an arc of intersection between the band joining $E$ to
itself and some $E_i$, an arc of intersection between the band and $E$ itself, or an arc of intersection of
the band with an annular collar.  Suppose first that the arc is not an arc where the band intersects
$E$ itself.  We verify that there is then a new boundary compression as in Case 1, 2, or 4, and we obtain a
contradiction in every case by returning to the analysis in each of those cases.  If the arc of intersection
is an arc where the band joining
$E$ to itself intersects
$E$, then the zipping cuts the band on the arc and rejoins it so that we get a shorter band joining $E$ to itself. 
The core of this sub-band is an arc $\gamma$ properly embedded in the new $\hat S=f\inverse(M)$.  In other words, if
the arc of intersection corresponds to an intersection point
$q$ in the interior of
$\beta'$ with
$\bdry E$, then
$q$ cuts
$\beta'$ into two arcs, including an arc $\gamma'$ joining $\bdry E$ to itself, which homotops to $\gamma$ in
$N(R_j)$.  If
$\gamma'\subset
\bdry R_j$ is not isotopic in
$\bdry R_j$ to an arc in
$\bdry E$, then clearly we get a new boundary compression along the the arc $\gamma$.  The
complexity
$n(\gamma')<n(\beta')$, so we replace $\beta'$ ($\beta$) by $\gamma'$ ($\gamma$).   If $\gamma'\subset \bdry R_j$ is
isotopic in
$\bdry R_j$ to an arc
$\delta'$ in
$\bdry E$, then the arc $\gamma$ in $\hat S$ is homotopic to an arc $\delta$ in the larger meridian disc $E$ in
$N(R_j)$.  Thus $\gamma\cup \delta$ is a simple closed curve in $\hat S$ which is null-homotopic in $M$, and by the
incompressibility of $\hat f$, the curve bounds a disc in $\hat S$.  It then follows that the original
$\beta'$ could be replaced by the closure in $\beta'$ of $\beta'-\gamma'$, which has smaller complexity and
corresponds to a boundary compression joining $E$ to itself.

Thus, we reduce to the case
that the zipping encounters no arc of intersection in an $R_j$. Therefore we are able to continue zipping
until the zipping branch point collides with another branch point, ends at a $K_i$ where
a collar meets $K_i$, or ends at $\bdry\bar M$.  The collision with another branch point
gives a surgery, and as usual, the surgery must cut a sphere from
$S$.  We obtain a new map of the same surface with fewer branch points, so the lexicographical complexity
$(c(f),b(f))$ is reduced.  In the special case where
$S$ is a sphere or a disc, we replace $f$ by the essential map on a sphere or disc component obtained from
surgery.  If a zipping
ends at a
$K_j$ as in Figures \ref{DehnZipCore} or \ref{DehnZipToK}, we obtain a contradiction to the minimal complexity of
$f$ as before, and similarly if the zipping ends at $\bdry \bar M$ as in Figure
\ref{DehnZipToBdry}.  

If at some point in the process, $P$ is the only remaining branch point, so that it is
impossible to start zipping at another branch point, then we begin zipping at $P$
and the zipping {\it must} end at a $K_i$ as in Figures \ref{DehnZipCore} or \ref{DehnZipToK}, or at $\bdry \bar M$
as in Figure \ref{DehnZipToBdry}. In all cases the 
complexity 
$c(f)$ is reduced, a contradiction.

Suppose now that we perform a sequence of zipping operations which cancel branch points
in pairs.  Ultimately, since the number of branch points is reduced after each cancellation of
branch points, a zipping must then end at
$P$ without first effecting any cut-and-paste's on arcs of $f(A)\cap (\cup_iE_i)$.  The
intersection of
$f(S)$ with
$R_j$ (Figure 13b) is unchanged by the zipping except on the arc ending at $P$.  The final
surgery eliminating
$P$ then does not separate
$S$.  This is a contradiction, since incompressibility guarantees that the surgery
corresponding to the cancellation of $P$ and another branch point cuts a sphere from $S$. 

\bigskip

\noindent {\it Case 6.} There remains the case in which a simple arc $\beta$ joins 
$\bdry_eS$ to itself and $\beta$ is homotopic in $M$ to an arc in $\bdry \bar M$.  Then by the
$\bdry$-incompressibility of $f$, $\beta$ cuts a half-disc $H$ from $S$.  Then $f(H)$ must intersect
$\cup K_j$, otherwise $\beta$ is not essential in $\hat S$, so replacing $f|_H$ by a map representing
the homotopy of the arc $\beta$ to $\bdry_eM$ reduces the complexity.

\bigskip

We must still consider the possibility that $\hat f\from \hat S\to M$ is a map of a disc or sphere.  We must
show that $\hat f$ is essential.  First suppose $\hat S$ is a 2-sphere, then $\hat S=S$ and $S$ is essential
in $\bar M$, so it is certainly essential in $M$.  Suppose that $\hat S$ is a disc with $\bdry \hat S\subset
\bdry_eM=\bdry \bar M$.  Then again $\hat S=S$. Since $S$ is essential in $\bar M$, $\hat S$ is essential in
$M$.  Next suppose $\hat S$ is a disc with $\bdry \hat S=\bdry_i\hat S$.  Then the essentiality of the map
$f\from (S,\bdry_iS)\to (\bar M, \cup_jK_j)$ implies the essentiality of $\hat f$.  Finally, suppose $\hat S$
is a disc with $\bdry \hat S=\bdry_m\hat S$.  Then $S$ is a sphere with $f\from S\to \bar M$ intersecting
$\cup_j K_j$ at a single point.  Thus $\hat f(\bdry\hat S)$ is a meridian in some $\bdry R_j$, and $\hat f$ is
essential.
\end{proof}

\section{Proofs of Other Results.}\label{Other}

\begin{proof} {\it (Proposition \ref{OldProp}).} Proposition \ref{OldProp} is a special case of Theorem
\ref{MainThm}, when
$\bdry
\bar M=\emptyset$ and $f\from S\to M$ is a map of a closed surface, so
$\bdry_eS=\bdry_i S=\emptyset$.
\end{proof}

\begin{proof} {\it (Proposition \ref{AgolProp})} Proposition \ref{AgolProp} is also a special case of Theorem
\ref{MainThm}.  Here
$M$ has a single torus boundary component, $\bar M=M(r)$, 
$\bdry M(r)=\emptyset$, and $S$ is a disc with $\bdry S=\bdry_iS$.  Since $f\from
(S,\bdry S) \to (M,K)$ is a map of a disc, where $K$ is the core curve of the Dehn
filling solid torus, relative incompressibility is automatic.
\end{proof}

The following is a slight variation of a result stating that any class in $H_2(M)$ can be represented by an
embedded oriented surface.

\begin{lemma} \label{HomologyToMapLemma} Let $Q$ be a 3-manifold, with $K_1,\ldots, K_n$ the components 
of an embedded (tame) link, $K=\cup_jK_j$.
Giving an orientation to each $K_j$, if $\ \sum_i n_{j}[K_{j}]=0$ in $H_1(Q)$ for some
integers
$n_{j}$, then there is a map of an oriented surface $f\from (S,\bdry S)\to (Q,K)$ with
$f_*([\bdry S])=\sum_j n_{j}[K_{j}]$ in $H_1(\coprod K_{j})$.  The map can be made an embedding
on the interior of $S$.  
\end{lemma}

\begin{proof}
The proof is similar to the proof of the result alluded to above, for
representing classes in homology by embedded surfaces.  In our situation we will not be
able to obtain an embedded surface, only a mapped-in surface.  We choose a triangulation
of $Q$ such that $K$ is a union of edges in the triangulation.  The fact that  $\sum_j
n_{j}[K_{j}]=0$ in $H_1(Q)$ gives a relative 2-cycle $z$ representing a class
in $H_2(Q,K)$, in simplicial homology.  Suppose
$z=\sum_i\sigma_i$, a finite sum, with
$\bdry z=\sum_j n_{j}[K_{j}]$.  We assume that the orientation of edges of $\partial z$ in $K_j$ coincides
with the orientation of $K_j$. We begin by pairing as many as possible cancelling pairs of
sides of 2-simplices in $z$ lying in
$K$, joining the pairs of 2-simplices on the interiors of the sides.  Our goal is not to obtain an embedded surface,
so the pairing can be done arbitrarily, subject to the condition that 2-simplices are joined with consistent
orientations.  For every other edge
$e$ of the triangulation, i.e., for every edge $e$ not contained in a $K_j$, we can pair
all occurrences of
$e$ in boundaries of 2-simplices, and join 2-simplices on the interiors of paired
sides, again joining 2-simplices so they have consistent orientations.  Thus we obtain a map
$g\from (P,\bdry P)\to (Q,K)$ of an oriented pseudo-surface $P$ with
$g(\bdry P)=\sum_j n_{j}[K_{j}]$, and with singularities mapped to the $0$-skeleton of
the triangulation.  Finally, making
$g$ transverse to 2-spheres linking the 0-simplices, the pre-image in $P$ is a collection
of simple closed curves and arcs with boundary in $\bdry P$. Cutting on these simple closed
curves and arcs detaches neighborhoods of singular points in $P$ to give a surface $\hat P$. 
Capping arcs and closed curves in
$\bdry \hat P$ mapped to linking spheres yields an oriented surface $S$ which can be mapped via $f\from
(S,\bdry S)\to (Q,K)$ to represent the same relative cycle as $g\from (P,\bdry P)\to (Q,K)$.  

To obtain a map which is an embedding on $\intr(S)$, one either chooses pairings more carefully away
from $K$, or one performs orientation-respecting cut-and-paste on curves of self-intersection away from $K$.
\end{proof}

Just as we can represent a class in $H_2(M)$ ($H_2(M,\bdry M)$) by an embedded
incompressible (and $\bdry$-incompressible) oriented surface, so also we can represent
a class in $H_2(Q,K)$ by an incompressible and relatively incompressible map of an
oriented surface, as we show in the following lemma.

\begin{lemma} \label{HomologyToEssentialLemma}  Let $Q$ be a 3-manifold, with $K_1,\ldots, K_n$ the components 
of an embedded (tame) link, $K=\cup_jK_j$.
Giving an orientation to each $K_j$, if $\sum_j n_{j}[K_{j}]=0$ in $H_1(Q)$ for some
integers
$n_{j}$, then there is a map of an oriented
surface
$f\from(S,\bdry S)\to (Q,K)$ with
$f_*([\bdry S])=\sum_j n_{j}[K_{j}]$, and with $f$ restricted to each component of 
$S$ essential.
\end{lemma}

\begin{proof}  By Lemma \ref{HomologyToMapLemma}, we can represent a class in $H_2(Q,K)$ whose boundary is
 $\sum_j n_{j}[K_{j}]$ by a map of an oriented surface $f\from (S,\bdry S)\to (Q,K)$ with
$f_*([\bdry S])=\sum_j n_{j}[K_{j}]$.  If $S$ has any closed components or null homologous components, we discard
these. Let us consider the restriction of
$f$ to a component
$F$ of
$S$ on which $f$ is not essential.  If the map is compressible, we can surger
$F$ to replace it with a surface $F'$ with $|\chi(F')|<|\chi(F)|$, or
if $F$ is an an annulus, we obtain a pair of discs $F'$.  At the same time we modify
$f$ to obtain
$f':(F',\bdry F')\to (Q,K)$ homologous to $f$.  If $f$ is relatively compressible, we can
surger
$F$ on a half-disc to obtain
$F'$ with
$|\chi(F')|<|\chi(F)|$ or if $F$ is an annulus, we obtain a disc $F'$. At the
same time we modify $f$ to obtain $f':(F',\bdry F')\to (Q,K)$ homologous to $f$. If
$F'$ has positive Euler characteristic, $f'$ is either inessential, hence null-homologous,
or it is essential. We discard the null-homologous components of $F'$.  We also discard closed components.  Finally
we replace $F$ by $F'$ to obtain a new $S$.  If $x(S)$ denotes the sum of $|\chi(G)|$ over components $G$ of $S$
having negative Euler characteristic, then $x(S)$ was reduced by our modification of $F$, or the number of inessential
components was decreased.   Repeating this process on any non-essential component of $S$, we must eventually obtain a
map which is  essential on each component of $S$.  
\end{proof}

\begin{proof}{\it(Corollary \ref{HomologyCor}.)}  We have the core curve $K$ of the Dehn filling solid torus
null-homologous in $M(r)$.  Hence, by Lemma  \ref{HomologyToMapLemma}, we can represent a generator in $H_2(M(r),
K)$ whose boundary is $[K]$ by a map $f\from (S,\bdry S)\to (M(r),K)$ of a surface $S$.
Actually, the proof of Lemma \ref{HomologyToMapLemma} shows that we can assume $\bdry S$ is mapped
homeomorphically to $K$.  By Lemma \ref{HomologyToEssentialLemma}, we can choose an essential map $f$ of a connected
surface with $f$ restricted to $\bdry S$ a homeomorphism to $K$.  The proof of Lemma 
\ref{HomologyToEssentialLemma} shows that the property of $f$ that it is a homeomorphism to $K$ on $\bdry S$ holds for
the new $f$.
  
Applying Theorem \ref{MainThm} gives a map $\hat f\from \hat S\to M$ as required, and $\hat S$ has the same genus as
$S$ by construction.
\end{proof}

\begin{proof}{\it (Corollary \ref{KnotCor})}  Corollary \ref{KnotCor} is similar to Corollary
\ref{HomologyCor}, but in this case we already know that there is a map $f\from S\to S^3$ with $\bdry f$ a
homeomorphism to $K$ and with $S$ a disc. Applying Theorem \ref{MainThm} immediately yields the
map $\hat f$ of a punctured disc, as described. 
\end{proof}

Corollary \ref{HomologyCorTwo} is a fairly obvious generalization of Corollary \ref{HomologyCor}.

\vfil
\eject
\bibliographystyle{amsplain}
\bibliography{ReferencesUO2}

\end{document}